\documentclass[12pt]{article}

\usepackage{a4wide}
\usepackage{amssymb}
\usepackage{amsfonts}
\usepackage{amsmath}

\date{}

\newtheorem{proposition}{Proposition}[section]
\newtheorem{theorem}[proposition]{Theorem}
\newtheorem{lemma}[proposition]{Lemma}

\newtheorem{corollary}[proposition]{Corollary}

\def\fd{ {\rm fd}}
\def\lfd{ {\rm lfd}}
\def\fdim{ {\rm d}}

\def\GK{{\rm  GK}\,}
\def\Kdim{{\rm K.dim }\,}

\def\Hom{{\rm Hom}}
\def\der{\partial }
\def\deri{{\partial }_i}
\def\derj{{\partial }_j}

\def\nFM0{{\nu }_{F,M_0}}
\def\nFN0{{\nu }_{F,N_0}}
\def\nGN0{{\nu }_{G,N_0}}

\def\nF{{\nu }_{F}}

\def\nG{ {\nu }_{G} }

\def\N0{ {\bf N}_0 }

\def\t{\otimes}
\def\g{\gamma}

\def\ra{\rightarrow}
\def\lra{\leftrightarrow}
\def\Xpm{X^{\pm }}

\def\l1{{\lambda}_1}

\def\a{\alpha}
\def\a0{ {\alpha }_0}
\def\a1{ {\alpha }_1}

\def\l{\lambda}
\def\o{\omega}

\def\nFGM0{{\nu }_{F,G,M_0}}

\def\lc{{\rm lc}}
\def\nFN0{{\nu}_{F,N_0}}

\def\lab{{\lambda }_{\alpha \beta }}
\def\Xa{X^{\alpha }}
\def\db{{\partial}^{\beta}}

\def\sm{{\sigma}^m}

\def\sm1{{\sigma}^{-1}}

\def\smtp1{{\sigma}^{-t+1}}

\def\o{\omega }
\def\S1{S^{-1}}

\def\Xpm1{X^{\pm 1}_1}

\def\sPM1{{\sigma }^{\pm 1}}
\def\sMP1{{\sigma }^{\mp 1 }}

\def\b{\beta}
\def\d{\delta}

\def\G{\Gamma}

\def\OO{{\cal O}}

\def\CD{{\cal D}}

\def\Ytm1{Y^{t-1}}
\def\Yim1{Y^{i-1}}

\def\CM{{\cal M}}

\def\CF{{\cal F}}

\def\Der{{\rm Der }}
\def\ad{{\rm ad }}
\def\dim{{\rm dim }}

\def\gr{ {\rm gr} }

\def\D{ \Delta }

\def\SL2Z{ {\rm SL}_2({\bf Z}) }

\def\Gp1{ G^{1 , 1 } }
\def\P11{ P^{-1 , 1 } }
\def\Pp1{ P^{1 , 1 } }

\def\nCLsr{{}^\nu\kern-2pt {\cal L}^{\sigma , \rho  }}
\def\nP{{}^\nu \kern-2pt P}
\def\nL{{}^\nu\kern-2pt L}
\def\nLL{{}^\nu\kern-2pt \Lambda}
\def\nPsr{{}^\nu\kern-2pt P^{\sigma , \rho  }}
\def\nLsr{{}^\nu\kern-2pt L^{\sigma , \rho  }}
\def\nuCL{{}^\nu\kern-2pt  {\cal L}}
\def\nCLsr{{}^\nu\kern-2pt {\cal L}^{\sigma , \rho  }}
\def\nCL1m{{}^\nu\kern-2pt {\cal L}^{-1 , 1  }}
\def\x1nu{x^\frac{1}{\nu}}
\def\xm1nu{x^{-\frac{1}{\nu}}}

\def\trdeg{{\rm tr.deg}}

\def\CB{{\cal B}}
\def\SL2Z{{\rm SL}_2( \mathbb{Z})}
\def\hSL2Z{\widehat{{\rm SL}}_2( \mathbb{Z})}

\def\bwnm1{ \overline{w}_n^{-1}}

\def\twnm1{ \widetilde{w}_n^{-1}}

\def\o0{\overline{0}}
\def\o1{\overline{1}}

\def\o{\omega}

\def\CB{{\cal B}}
\def\SL2Z{{\rm SL}_2( \mathbb{Z})}
\def\hSL2Z{\widehat{{\rm SL}}_2( \mathbb{Z})}

\def\bwnm1{ \overline{w}_n^{-1}}

\def\twnm1{ \widetilde{w}_n^{-1}}

\def\b0{ \overline{0}}

\def\Gp{\mathfrak{p}}

\def\hol{{\rm hol}}
\def\Cdim{{\rm Cdim}}
\def\CDX{{\cal D}(X)}
\def\GKtrdeg{{\rm GKtr.deg}}

\def\CD{{\cal D}}
\def\End{{\rm End}}
\def\D{{\Delta}}
\def\gr{{\rm gr}}

\def\a{\alpha}
\def\b{\beta}

\def\CD{{\cal D}}

\begin{document}

\author{V.\  Bavula  }

\title{Filter Dimension}

\maketitle
{\bf Contents}
\begin{enumerate}
\item Introduction. \item Filter dimension of algebras and
modules. \item The first  filter inequality. \item Krull,
Gelfand-Kirillov and filter dimensions of simple finitely
generated algebras.\item Filter dimension of the ring of
differential operators on a smooth irreducible affine algebraic
variety.  \item The multiplicity for the filter dimension,
holonomic modules over simple finitely generated algebras. \item
Filter dimension and commutative subalgebras of simple finitely
generated algebras and their division algebras. \item Filter
dimension and isotropic subalgebras of Poisson algebras.
\end{enumerate}


\section{Introduction}
Throughout the paper, $K$ is a  field, a module $M$ over an
 algebra $A$  means a {\em left} module denoted ${}_AM$, $\t
=\t_K$.

 Intuitively, the filter dimension of an algebra or a
module measures how `close' standard filtrations of the algebra or
the module are. In particular, for a simple algebra it also
measures the growth of how `fast' one can prove that the algebra
is simple.

The filter dimension appears naturally when one wants to
generalize the Bernstein's inequality for the Weyl algebras to the
class of simple finitely generated algebras.

   The $n$'th {\bf Weyl algebra} $A_n$ over the field $K$ has
    $2n$ generators $X_1,\ldots ,X_n$, ${\der}_1,\ldots ,{\der
}_n$ that satisfy the  defining  relations
$$\deri X_j-X_j\deri =\delta {}_{ij},\,\hbox{the Kronecker delta},
\,X_iX_j-X_jX_i=\deri \derj -\derj \deri=0,$$ for all
$i,j=1,\ldots ,n$. When  char $K=0$ the Weyl algebra $A_n$ is a
simple Noetherian
 finitely generated algebra canonically isomorphic to the ring of
differential operators $K[X_1,\ldots ,X_n,\frac{\der}{\der
X_1},\ldots ,\frac{\der }{\der X_n}]$ with polynomial coefficients
($X_i\lra X_i, {\der}_i\lra \frac{\der_i}{\der X_i}$, $i=1,\ldots
,n$).

Let $\Kdim $ and $\GK $ be the {\bf (left) Krull}  (in the sense
of Rentschler and Gabriel, \cite{Ren-Gab}) and the {\bf
Gelfand-Kirillov} dimension respectively.

\begin{theorem}\label{Berin}
({\bf The Bernstein's  inequality}, \cite{Ber72}) Let $A_n$ be the
$n$'th Weyl algebra over a field of characteristic zero. Then $\GK
(M)\geq n$ for all nonzero finitely generated $A_n$-modules $M$.
\end{theorem}

Let $A$ be a simple finitely generated infinite dimensional
$K$-algebra. Then $\dim_K(M)=\infty $ for all nonzero $A$-modules
$M$ (the algebra $A$ is simple, so the $K$-linear map $A\ra {\rm
Hom}_K(M,M)$, $a \mapsto (m \mapsto am)$, is injective, and so
$\infty =\dim_K(A)\leq  \dim_K({\rm Hom}_K(M,M))$ hence
$\dim_K(M)=\infty$). So, the Gelfand-Kirillov dimension (over $K$)
$\GK (M)\geq 1$ for all nonzero $A$-modules $M$.

{\it Definition}. $h_A:=\inf \{ \GK (M)\, | \,  M $  is a nonzero
finitely  generated $A$-module$\}$ is called the {\bf holonomic
number} for the algebra $A$.

{\sc Problem}. {\em For a simple finitely generated algebra find
its holonomic number.}

 To find  an
approximation of the holonomic number for simple finitely
generated  algebras and to generalize the Bernstein's inequality
for these algebras was a main motivation for introducing the
filter dimension, \cite{Bavcafd}. In this paper $d$ stands for the
filter dimension $\fd$ or the left filter dimension $\lfd$. The
following two inequalities are central for the proofs of almost
all results in this paper.

\begin{itemize}
\item {\bf The First Filter Inequality, \cite{Bavcafd}.}  {\em Let
$A$ be a simple finitely generated algebra. Then
$$
 \GK (M)\geq \frac{\GK (A)}{\fdim (A)+\max \{ \fdim (A),
1\} }
$$
 for all nonzero finitely generated $A$-modules $M$.}

\item {\bf The Second Filter Inequality, \cite{Bavjafd}.}  {\em
Under a certain mild conditions (Theorem \ref{SFI}) the (left)
Krull dimension of the algebra $A$ satisfies the following
inequality}
$$
{\rm K.dim} (A)\leq \GK (A)(1-\frac{1}{\fdim (A)+\max \{ \fdim
(A), 1\}}).
$$

\end{itemize}

The paper is organized as follows. Both filter dimensions are
introduced in Section \ref{fdamhaz}. In Sections \ref{IandIIhaz}
 and \ref{KGKFdsafhaz} the first and the second filter inequalities
 are proved respectively. In Section \ref{KGKFdsafhaz} we use both
 filter inequalities for giving short proofs of some classical
 results about the rings $\CDX$ of differential operators on
 smooth irreducible affine algebraic varieties. The (left) filter
 dimension of $\CDX $ is $1$ (Section \ref{5fdiffhaz}). A concept
 of multiplicity for the filter dimension and a concept of
 holonomic module for (simple) finitely generated algebras appear
 in Section \ref{multholhaz}. Every holonomic module has finite
 length (Theorem \ref{holfhaz}). In Section \ref{hazmaxissPA}
 an upper bound is given $(i)$ for the
Gelfand-Kirillov dimension of commutative subalgebras of simple
finitely generated infinite dimensional algebras (Theorem
\ref{GKcsuba}), and $(ii)$ for the transcendence degree of
 subfields of quotient  rings of (certain) simple finitely
generated infinite dimensional algebras (Theorems \ref{trdsubf}
and \ref{c1trds}). In Section \ref{hazisotr} a similar upper bound
is obtained for the Gelfand-Kirillov dimension of {\em isotropic}
subalgebras of strongly simple  Poisson algebras (Theorem
\ref{PoGKcsu}).


\section{Filter dimension of algebras and
modules}\label{fdamhaz}

In this section, the filter dimension of algebras and modules will
be defined.

{\bf The Gelfand-Kirillov dimension}. Let $\CF $ be the set of all
functions from the set of natural numbers $ \mathbb{N}=\{ 0, 1,
\ldots \}$ to itself. For each function $f\in \CF $, the
non-negative real number or $\infty $ defined as
$$ \g (f):=\inf \{  r\in \mathbb{R}\, | \, f(i)\leq i^r\; {\rm for
}\; i\gg 0\}$$ is called the  {\bf degree} of $f$. The function
$f$ has {\bf polynomial growth} if $\g (f)<\infty $. Let $f,g,
p\in \CF $, and $p(i)=p^*(i)$ for $i\gg 0$ where $p^*(t)\in
\mathbb{Q}[t]$ (a polynomial algebra with coefficients from the
field  of rational numbers). Then
\begin{eqnarray*}
\g (f+g)\leq \max \{ \g  (f), \g (g)\}, & & \g (fg)\leq \g (f)+ \g
(g),\\
\g (p)=\deg_t(p^*(t)), & & \g (pg)= \g (p)+ \g
(g).\\
\end{eqnarray*}
Let $A=K\langle a_1, \ldots , a_s\rangle $ be a finitely generated
 $K$-algebra. The finite dimensional filtration $F=\{ A_i\}$
associated with algebra generators $ a_1, \ldots , a_s$:
$$ A_0:=K\subseteq A_1:=K+\sum_{i=1}^sKa_i\subseteq \cdots \subseteq
A_i:=A_1^i\subseteq \cdots $$ is  called the {\bf standard
filtration} for the algebra $A$. Let $M=AM_0$ be a finitely
generated $A$-module where $M_0$ is a finite dimensional
generating subspace. The finite dimensional filtration $\{
M_i:=A_iM_0\}$ is called the {\bf standard filtration} for the
$A$-module $M$.

{\it Definition}. $\GK (A):=\g (i\mapsto \dim_K(A_i))$ and $ \GK
(M):=\g (i\mapsto \dim_K(M_i))$ are called the {\bf
Gelfand-Kirillov} dimensions of the algebra $A$ and the $A$-module
$M$ respectively.

It is easy to prove that the Gelfand-Kirillov dimension  of the
algebra (resp. the module)  does not depend on the choice of the
standard filtration of the algebra (resp. and the choice of the
generating subspace of the module).

{\bf The return functions and the (left) filter dimension}.

 {\it Definition} \cite{Bavcafd}. The function $\nFM0 :\mathbb{N}\ra
\mathbb{N}\cup \{ \infty \}$,
$$ \nFM0 (i):=\min \{ j\in\mathbb{N}\cup \{ \infty \}: \;
A_jM_{i,gen}\supseteq M_0\;\; {\rm for \; all}\;\; M_{i,gen} \}$$
is called the {\bf return function} of the $A$-module $M$
associated with the filtration $F=\{ A_i\}$ of the algebra $A$ and
the generating subspace $M_0$ of the $A$-module $M$ where
$M_{i,gen}$ runs through all generating subspaces for the
$A$-module $M$ such that $M_{i,gen}\subseteq M_i$.

Suppose, in addition, that the finitely generated algebra $A$ is a
{\em simple} algebra.   The {\bf return function} $\nu_F \in \CF $
and the {\bf left return function} $\l_F\in \CF $ for the algebra
$A$ with respect to the standard filtration $F:= \{ A_i\}$ for the
algebra $A$ are defined by the rules:
\begin{eqnarray*}
\nu_F(i)&:=& \min \{ j\in \mathbb{N}\cup \{ \infty \} \, | \,\,
1\in A_jaA_j\;\, {\rm
for \; all}\; \,  0\neq a\in A_i\},\\
\l_F(i)&:=& \min \{ j\in \mathbb{N}\cup \{ \infty \} \, | \, 1\in
AaA_j\; \; {\rm for \; all}\; \; 0\neq a\in A_i\},
\end{eqnarray*}
where $A_jaA_j$ is the vector subspace of the algebra $A$ spanned
over the field $K$ by the elements $xay$ for all $x,y\in A_j$; and
$AaA_j$ is the left ideal of the algebra $A$ generated by the set
$aA_j$.  The next result shows that under a mild restriction the
maps $\nu_F(i)$ and $\l_F(i)$ are finite.

Recall that the centre of a simple algebra is a field.
\begin{lemma}\label{fdfin}
Let $A$ be a simple finitely generated algebra such that its
centre $Z(A)$ is an algebraic field extension of $K$. Then
 $\l_F(i)\leq \nu_F(i)<\infty $ for all $i\geq 0$.
\end{lemma}

{\it Proof}. The first inequality is evident.

 The centre $Z=Z(A)$ of the simple algebra $A$ is a field that
 contains $K$. Let $\{ \o_j\, | \, j\in J\}$ be a $K$-basis for the $K$-vector
space $Z$. Since $\dim_K(A_i)<\infty $, one can find  finitely
many $Z$-linearly independent elements, say $a_1,\ldots , a_s$, of
$A_i$ such that $A_i\subseteq Za_1+\cdots +Za_s$. Next, one can
find a finite subset, say $J'$, of $J$ such that $A_i\subseteq
Va_1+\cdots +Va_s$ where $V=\sum_{j\in J'}K\o_j$. The field $K'$
generated over $K$ by the elements $\o_j$, $j\in J'$, is a finite
field extension of $K$ (i.e. $\dim_K(K')<\infty $) since $Z/K$ is
algebraic, hence $K'\subseteq A_n$ for some $n\geq 0$. Clearly,
$A_i\subseteq K'a_1+\cdots +K'a_s$.

The $A$-bimodule ${}_AA_A$ is simple with ring of endomorphisms
${\rm End}({}_AA_A)\simeq Z$. By the {\em Density} Theorem,
\cite{Pierceb}, 12.2, for each integer $1\leq j \leq s$, there
exist elements of the algebra $A$, say $x_1^j, \ldots , x_m^j,
y_1^j, \ldots , y_m^j$, $m=m(j)$, such that for all $1\leq l\leq
s$
$$ \sum_{k=1}^m x_k^ja_ly_k^j=\delta_{j,l}, \;\; {\rm the \; Kronecker
\; delta}.$$ Let us fix a natural number, say $d=d_i$, such that
$A_d$ contains all the elements $x_k^j$, $y_k^j$, and the field
$K'$. We claim that $\nu_F(i)\leq 2d$. Let $0\neq a\in A_i$. Then
$a=\l_1a_1+\cdots +\l_sa_s$ for some $\l_i\in K'$. There exists
$\l_j\neq 0$. Then $\sum_{k=1}^m\l_j^{-1}x_k^ja_jy^j_k=1$, and $
 \l_j^{-1}x_k^j, y^j_k\in A_{2d}$. This proves the claim and the lemma.
  $\Box $

{\it Remark}.  If the field $K$ is {\em uncountable} then
automatically the centre $Z(A)$ of a simple finitely generated
algebra $A$ is algebraic over $K$  (since $A$ has a countable
$K$-basis and the rational function field $K(x)$ has uncountable
basis over $K$ since elements $\frac{1}{x+\l }$, $\l \in K$, are
$K$-linearly independent).

It is easy to see that for a finitely generated  algebra $A$ any
two standard finite dimensional filtrations $F=\{ A_i\} $ and
$G=\{ B_i \}$ are {\bf equivalent}, $(F\sim G),$ that  is, there
exist natural numbers $a,b,c,d$ such that
$$A_i\subseteq B_{ai+b}\;\hbox{and}\;B_i\subseteq A_{ci+d}\;\hbox{for}\;
i\gg 0 .$$

If one of the inclusions holds, say the first, we write $F\leq G$.

\begin{lemma}      
 Let $A$ be a finitely generated  algebra equipped with two standard
finite dimensional filtrations $F=\{ A_i\}$ and $G =\{ B_i\}.$
\begin{enumerate}
\item  Let $M$ be a finitely generated $A$-module. Then $\g (\nFM0
)=\g (\nGN0 )$ for any finite dimensional generating subspaces
$M_0$ and $N_0$ of the $A$-module $M$. \item  If, in addition, $A$
is a simple algebra  then  $\g (\nF )=\g (\nG )$ and $\g (\l_F)=\g
(\l_G)$.
\end{enumerate}
\end{lemma}

{\it Proof}. 1. The module $M$ has two standard finite dimensional
filtrations $\{M_i=A_iM_0 \}$ and $\{N_i=B_iN_0 \} $. Let $\nu
=\nFM0 $    and $\mu =\nGN0 $.

Suppose that $F=G$. Choose a natural number $s$ such that
$M_0\subseteq N_s$ and $N_0\subseteq M_s$, so $N_i\subseteq
M_{i+s}$ and $M_i\subseteq N_{i+s}$ for all $i\geq 0$. Let
$N_{i,gen}$ be any generating subspace for the $A$-module $M$ such
that $N_{i,gen}\subseteq N_i$.  Since $M_0\subseteq A_{\nu
(i+s)}N_{i,gen}$ for all $i\geq 0$ and $N_0 \subseteq A_sM_0,$ we
have $N_0\subseteq A_{\nu (i+s)+s}N_{i,gen},$ hence, $\mu (i)\leq
\nu (i+s)+s$ and finally $\g (\mu )\leq \g (\nu ).$ By symmetry,
the opposite inequality is true and  so $\g (\mu )=\g (\nu ) $.

Suppose that $M_0=N_0$. The  algebra $A$ is  a finitely generated
algebra, so all standard finite dimensional filtrations of the
algebra $A$  are equivalent. In particular,  $F\sim G$  and so one
can choose natural numbers $a,b,c,d$ such that
$$A_i\subseteq B_{ai+b}\;\hbox{and}\;B_i\subseteq A_{ci+d}\;\hbox{for}\;
i\gg 0 .$$

 Then $N_i=B_iN_0\subseteq A_{ci+d}M_0=M_{ci+d}$ for all $i\geq
 0$, hence $N_0=M_0\subseteq A_{\nu (ci+d)}N_{i,gen}\subseteq B_{a\nu
(ci+d)+b} N_{i,gen}$, therefore $\mu (i)\leq a\nu (ci+d)+b$ for
all $i\geq 0$,  hence $\g (\mu )\leq \g (\nu ) $. By symmetry, we
get
 the  opposite inequality which implies  $\g (\mu )=\g (\nu )$.
 Now, $\g (\nFM0 )=\g (\nFN0 )=\g (\nGN0 ) $.

2. The algebra $A$ is simple, equivalently, it is  a simple (left)
$A\t A^0$-module where $A^0$  is the opposite algebra to $A$. The
opposite algebra has the standard filtration $F^0=\{ A_i^0\}$,
opposite to the filtration $F$. The tensor product of algebras
$A\t A^0$, so-called,  the enveloping algebra of $A$, has the
standard filtration $F\t F^0= \{ C_n\}$ which is the tensor
product of the standard filtrations $F$ and $F^0$, that is,
$C_n=\sum \{ A_i\t A_j^0,i+j\leq n \} $. Let  ${\nu }_{ F\t F^0, K
}$ be the return function of the $A\t A^0$-module $A$ associated
with the filtration $F\t F^0$ and the generating subspace $K$.
Then
$$\nF (i)\leq {\nu }_{ F\t F^0,K }(i)\leq 2\nF (i) \; \hbox{for all}\;  i\geq 0,$$
 and so
\begin{equation}\label{gnuFK}
\g (\nF )=\g ({\nu }_{F\t F^0,K}),
\end{equation}
 and, by the first statement, we
have  $\g (\nu_F)= \g ({\nu }_{F\t F^0,K})=\g ({\nu }_{G\t
G^0,K})=\g (\nu_G)$, as required. Using a similar argument as in
the proof of the first statement one can proof that $\g (\l_F)=\g
(\l_G)$. We leave this as an exercise. $\Box $

{\it Definition} \cite{Bavcafd}. $\fd (M)=\g (\nFM0 )$ is the {\bf
filter dimension} of the $A$-module $M$, and $\fd (A):= \fd
({}_{A\t A^0}A)$ is the {\bf filter dimension} of the algebra $A$.
If, in addition, the algebra $A$ is simple, then $\fd (A)=\g
(\nu_F)$, and $\lfd (A):=\g (\l_F)$ is called the {\bf left filter
dimension} of the algebra $A$.

By the previous lemma the definitions make sense  (both filter
dimensions do not depend on the choice of the standard filtration
$F$ for the algebra $A$).

By Lemma \ref{fdfin}, $\lfd (A)\leq \fd (A)$.

{\it Question. What is the filter dimension of a polynomial
algebra?}


\section{The first filter inequality}\label{IandIIhaz}

 In this paper,  $\fdim (A)$
means either the filter dimension $\fd (A)$ or the left filter
dimension $\lfd (A)$ of a simple finitely generated  algebra $A$
(i.e. $\fdim =\fd , \lfd $).
  Both filter dimensions  appear
  naturally when one tries to find a
{\em lower} bound for the holonomic number (Theorem \ref{FFI}) and
an {\em upper} bound (Theorem \ref{SFI}) for the (left and right)
{\em Krull} dimension (in the sense of Rentschler-Gabriel,
\cite{Ren-Gab}) of simple finitely generated algebras.

The next theorem is a generalization of the {\bf Bernstein's
Inequality} (Theorem \ref{Berin}) to the  class of simple finitely
generated algebras.

\begin{theorem}\label{FFI}

({\bf The First Filter Inequality},  \cite{Bavcafd, bie98})
 Let $A$ be a simple finitely generated  algebra.
 Then
$$
 \GK (M)\geq \frac{\GK (A)}{\fdim(A)+\max \{ \fdim(A),
1\} }
$$
 for all nonzero finitely generated $A$-modules $M$ where $\fdim
 =\fd , \lfd$.
\end{theorem}

{\it Proof}. Let $\l =\l_F$ be the left return function associated
with a standard filtration $F$ of the algebra $A$ and let  $0\neq
a \in A_i$. It suffices to prove the inequality for $\l $ (since
$\fd (A)\geq \lfd (A)$). It follows from the inclusion
$$AaM_{\l (i)}=AaA_{\l (i)}M_0\supseteq 1M_0=M_0$$
that the linear map
$$A_i\ra  \Hom (M_{\l (i)},M_{\l (i)+i}), a\mapsto
(m\mapsto  am),$$ is injective, so  dim $A_i\,\le $ dim $M_{\l
(i)}$ dim $M_{\l  (i)+i}$.
 Using the above elementary properties of the degree  (see also
\cite{MR}, 8.1.7), we have
\begin{eqnarray*}
\GK (A)& =& \gamma (\dim \, A_i)\le \gamma (\dim \, M_{\l
(i)})+\gamma
(\dim \, M_{\l (i)+i})\\
&\le &  \gamma (\dim \, M_i )\gamma (\l )+\gamma (\dim \, M_i
)\max \{
\gamma (\l ),1\}\\
&=& \GK (M)(\lfd A+\max \{\lfd A, 1\})\\
&\leq & \GK (M)(\lfd A+\max \{\lfd A, 1\}). \;\; \Box
\end{eqnarray*}

The result above gives a lower bound for the holonomic number of a
simple finitely generated algebra
$$ h_A\geq  \frac{\GK (A)}{\fdim(A)+\max \{ \fdim(A),
1\} }.
$$
\begin{theorem}\label{IIIin}
Let $A$ be a finitely generated  algebra. Then
$$\GK (M)\leq \GK(A)\,\fd (M)$$
for any simple $A$-module $M$.
\end{theorem}

{\it Proof}. Let $\nu =\nu_{F, Km} $ be the return function of the
module  $M$ associated with a standard finite dimensional
filtration $F=\{A_i\}$ of the algebra $A$ and a fixed nonzero
element  $m\in M$. Let $\pi :M\ra K$ be a non-zero linear map
satisfying $\pi (m)=1$. Then, for any $i\geq 0$ and any $0\neq
u\in M_i$: $1=\pi (m)\in \pi (A_{\nu (i)}u)$, and so the linear
map
$$M_i\ra \Hom (A_{\nu (i)},K),\,\,u\mapsto (a\mapsto \pi (au)),$$
is an injective map hence dim $M_i\leq $ dim $A_{\nu (i)}$ and
finally  $\GK (M)\leq \GK (A)\, \fd (M)$.  $\Box$

\begin{corollary}
Let $A$ be a simple finitely generated  infinite dimensional
algebra. Then
 $$\fd (A)\geq \frac{1}{2}.$$
\end{corollary}

{\it Proof}. The algebra $A$ is a finitely generated infinite
dimensional algebra hence $\GK (A)>0$. Clearly,
 $\GK(A\t A^0)\leq \GK (A)+\GK (A^0)=2\GK(A)$.
Applying Theorem  \ref{IIIin} to  the  simple $A\t A^0$-module
$M=A$ we finish the proof:
$$ \GK (A)=\GK ({}_{A\t A^0}A)\leq \GK (A\t A^0)\fd
({}_{A\t A^0}A)\leq 2\GK(A) \fd (A)$$ hence $\fd (A)\geq
\frac{1}{2}$. $\Box $

{\it Question. Is $\fd (A)\geq 1$ for all simple finitely
generated infinite dimensional algebras $A$?}

{\it Question. For which numbers $d\geq \frac{1}{2}$ there exists
a  simple finitely generated infinite dimensional algebra $A$ with
$\fd (A)=d$?}

\begin{corollary}
Let $A$ be a simple finitely generated  infinite dimensional
algebra. Then
 $$\fd (M)\geq \frac{1}{\fd (A)+\max \{ \fd (A), 1\} }$$
  for all simple $A$-modules $M$.
\end{corollary}

{\it Proof}. Applying Theorem \ref{FFI} and Theorem \ref{IIIin},
we have the result
$$ \fd (M)\geq \frac{\GK (M)}{\GK (A)}\geq
\frac{\GK (A)}{\GK (A) (\fd (A)+\max \{ \fd (A), 1\}
)}=\frac{1}{\fd (A)+\max \{ \fd (A), 1\} }. \;\; \Box
$$


\section{Krull, Gelfand-Kirillov and filter dimensions of simple
finitely generated
algebras}\label{KGKFdsafhaz}

In this section,  we prove the second filter inequality (Theorem
\ref{SFI}) and apply both filter inequalities for giving short
proofs of some classical results about the rings of differential
operators on a smooth irreducible affine algebraic varieties
(Theorems \ref{Berin}, \ref{KdimDB}, \ref{BerinDB},
\ref{KdimAnRGab}).

We say that an algebra $A$ is {\bf (left) finitely partitive}
 (\cite{MR}, 8.3.17)
if, given any finitely generated $A$-module $M$, there is an
integer $n=n(M)>0$ such that for every strictly descending chain
of $A$-submodules of $M$:
$$M=M_0\supset M_1\supset \cdots \supset M_m$$
with $\GK (M_i/M_{i+1})=\GK (M)$, one has $m\leq n$. McConnell and
Robson write in their book \cite{MR}, 8.3.17,  that ``{\em yet  no
examples are known which fail  to have this property.}''

Recall that $\Kdim $ denotes the (left) Krull dimension in the
sense of Rentschler and Gabriel, \cite{Ren-Gab}.

\begin{lemma}\label{Mab}
Let $A$ be a finitely partitive algebra with $\GK (A)<\infty $.
 Let $a\in  \mathbb{N}$, $b\geq 0$ and suppose that
$\GK (M)\geq a+b$ for all finitely generated  $A$-modules $M$ with
${\rm K.dim}(M)=a$,   and  that $\GK (N)\in \mathbb{N}$ for all
finitely generated $A$-modules $N$ with ${\rm K.dim}(N)\geq a$.
Then $\GK (M)\geq {\rm K.dim} (M)+b$ for all finitely generated
$A$-modules $M$ with $\Kdim (M)\geq a$. In particular, $\GK
(A)\geq {\rm K.dim} (A)+b$.
\end{lemma}

{\it Remark.} It is assumed that a module $M$ with $\Kdim (M)=a$
exists.

 {\it Proof}. We use induction on $n={\rm K.dim}(M)$. The base
of induction, $n=a$, is true. Let $n>a$. There exists a descending
chain of submodules $ M=M_1\supset M_2\supset \cdots $ with
 ${\rm K.dim}(M_i/M_{i+1})=n-1$ for $i\geq 1$. By induction,
 $\GK (M_i/M_{i+1})\geq n-1+b$ for  $i\geq 1$. The algebra $A$ is
finitely partitive, so  there exists $i$ such that $\GK (M)>\GK
(M_i/M_{i+1})$, so
 $\GK (M)-1\geq \GK (M_i/M_{i+1})\geq n-1+b$, since
  $\GK (M)\in  \mathbb{N}$, hence $\GK (M)\geq {\rm K.dim} (M)+b$.
  Since $\Kdim (A)\geq \Kdim (M)$ for all finitely generated
  $A$-modules $M$ we have $\GK (A)\geq {\rm K.dim} (A)+b$.
   $\Box $

\begin{theorem}\label{SFI}
(\cite{Bavjafd})  Let $A$ be a simple   finitely generated
finitely partitive algebra with $\GK (A)<\infty $. Suppose that
the Gelfand-Kirillov dimension of every finitely generated
$A$-module is a natural number.
 Then
$$
{\rm K.dim}(M)\leq \GK (M) - \frac{\GK (A)}{{\rm d}(A)+\max \{{\rm
d}(A), 1\}}
$$
for any nonzero finitely generated $A$-module $M$.  In particular,
$$
{\rm K.dim} (A)\leq \GK (A)(1-\frac{1}{{\rm d}(A)+\max \{{\rm
d}(A), 1\}}).
$$
\end{theorem}

{\it Proof}.  Applying the lemma above to the family of finitely
generated $A$-modules of Krull dimension $0$, by Theorem
\ref{FFI}, we can put $a=0$ and
$$b=\frac{\GK (A)}{{\rm d}(A)+\max
\{{\rm d}(A), 1\}},$$ and the result follows. $\Box $


Let $K$ be a field of characteristic zero and $B$ be a commutative
$K$-algebra. The ring of ($K$-linear) {\bf differential operators}
$\CD (B)$ on $B$ is defined as $\CD (B)=\cup_{i=0}^\infty \,\CD_i
(B)$ where $\CD_0 (B)=\End_R(B)\simeq B$, ($(x\mapsto bx)\lra b$),
$$ \CD_i (B)=\{ u\in \End_K(B):\, [u,r]\in \CD_{i-1} (B)\; {\rm for \; each \; }\; r\in B\}.$$
Note that the $\{ \CD_i (B)\}$ is, so-called, the {\bf order
filtration} for the algebra $\CD (B)$:
$$\CD_0(B)\subseteq   \CD_1 (B)\subseteq \cdots \subseteq
\CD_i (B)\subseteq \cdots\;\; {\rm and}\;\; \CD_i (B)\CD_j
(B)\subseteq \CD_{i+j} (B), \;\; i,j\geq 0.$$

The subalgebra $\D (B)$ of $\End_K(B)$ generated by $B\equiv
\End_R(B)$ and by the set ${\rm Der}_K (B)$ of all $K$-derivations
of $B$ is called the {\bf derivation ring} of $B$. The derivation
ring $\D (B)$ is a subring of  $\CD (B)$.

Let the {\em finitely generated} algebra $B$ be a {\em regular
commutative domain of Krull dimension} $n<\infty $. In geometric
terms, $B$ is the coordinate ring $\OO (X)$ of a smooth
irreducible  affine algebraic variety $X$ of dimension $n$. Then
\begin{itemize}
\item ${\rm Der}_K (B)$ {\em is a finitely generated projective}
$B$-{\em module of rank} $n$; \item  $\CD (B)=\Delta (B) $; \item
$\CD (B)$ {\em is a simple (left and right) Noetherian domain with
}  $\GK \, \CD (B)=2n$ ($n=\GK (B)=\Kdim (B))$; \item  $\CD
(B)=\Delta (B)$ {\em is an almost centralizing extension of} $B$;
\item {\em the associated graded ring} $\gr \, \CD (B) =\oplus \,
\CD_i(B)/\CD_{i-1}(B)$ {\em is a commutative domain}; \item {\em
the Gelfand-Kirillov dimension of every finitely generated}
 $\CD (B) $-{\em module is a natural number}.
\end{itemize}

For the proofs of the statements above the reader is referred to
\cite{MR}, Chapter 15.
 So, the domain $\CD (B)$ is a simple finitely generated infinite dimensional Noetherian algebra
(\cite{MR}, Chapter 15).

{\em Example}. Let $P_n=K[X_1, \ldots ,X_n]$ be a polynomial
algebra.
 ${\rm Der}_K (P_n)=\oplus_{i=1}^n \; P_n \frac{\partial }{\partial X_i}$,
$$ \CD (P_n)=\Delta (P_n)  = K[X_1, \ldots ,X_n,  \frac{\partial }{\partial X_1}, \ldots ,
\frac{\partial }{\partial X_n} ]$$ is the ring of differential
operators with polynomial coefficients, i.e.    the $n$'th Weyl
algebra $A_n$.

In Section 5, we prove the following result.

\begin{theorem}\label{fdif=1}
(\cite{Bavjafd}) The  filter dimension  and the left filter
dimension  of the ring of differential operators $\CD (B)$ are
both  equal to $1$.
\end{theorem}

As an application we compute the Krull dimension of $\CD (B) $.

\begin{theorem}\label{KdimDB}
(\cite{MR}, Ch. 15)
$$
{\rm K.dim} \, \CD (B)=\frac{\GK (\CD (B))}{2} ={\rm K.dim} (B).
\;
$$
\end{theorem}

{\it Proof}. The second equality is clear
 $(\GK (\CD (B))=2 \GK (B)=2{\rm K.dim} (B))$.
 It follows from Theorems \ref{SFI} and \ref{fdif=1} that
 $${\rm K.dim} \, \CD (B)\leq \frac{\GK (\CD (B))}{2} =\Kdim (B). $$ The map $I\ra \CD (B)I$ from
 the set of  left ideals of $B$ to the set of  left ideals  of $\CD (B)$ is injective, thus
 ${\rm K.dim} \, (B)\leq {\rm K.dim}\, \CD (B)$. $\Box $

This result shows that for the ring of differential operators on a
smooth irreducible affine algebraic  variety the inequality in
Theorem \ref{SFI} is {\em the equality}.

\begin{theorem}\label{BerinDB}
 (\cite{MR}, 15.4.3) Let $M$ be a nonzero finitely generated $\CD (B)$-module. Then
$$ \GK (M)\geq \frac{\GK (\CD (B))}{2} ={\rm K.dim} (B).  $$
\end{theorem}

{\it Proof}. By Theorems \ref{FFI} and \ref{fdif=1},
$$ \GK (M)\geq \frac{\GK (\CD (B))}{1+1}=\frac{2\GK (B)}{2}=\GK
(B)=\Kdim (B).\;\;  \Box $$

So, for the ring of differential operators on a smooth affine
algebraic variety the inequality in Theorem \ref{FFI} is in fact
{\em the equality}.

In general, it is difficult to find the exact value for the filter
dimension but for the Weyl algebra $A_n$ it is easy and one can
find it directly.

\begin{theorem}\label{fdWeyl=1}
Both the  filter dimension  and the left filter dimension  of the
Weyl algebra $A_n$ over a field of characteristic zero  are equal
to $1$.
\end{theorem}

{\it Proof}. Denote by $a_1, \ldots ,a_{2n}$ the canonical
generators of the Weyl algebra $A_n$ and denote by $F=\{
A_{n,i}\}_{i\geq  0}$ the standard filtration associated with the
canonical generators. The associated graded algebra $\gr \,
A_n:=\oplus_{i\geq 0}\,A_{n,i}/A_{n,i-1}$, $( A_{n,-1}=0)$ is a
polynomial algebra in $2n$ variables, so
$$
\GK (A_n)=\GK (\gr \,A_n)=2n.
$$
For every $i\geq 0$:
$$ \ad \,a_j:A_{n,i}\ra  A_{n,i-1}, \; x\mapsto
\ad \,a_j(x):=a_jx-xa_j.
$$
The algebra $A_n$ is central ($Z(A_n)=K$), so
$$
\ad \,a_j(x)=0 \; {\rm for \; all}\; j=1,\ldots ,2n
\Leftrightarrow
 x\in Z(A_n)=K=A_{n,0}.
$$
These two facts imply $\nu_F(i)\leq  i$ for $i\geq 0$, and  so
 ${\rm d }(A_n)\leq 1$.

The $A_n$-module $P_n:=K[X_1, \ldots , X_n]\simeq
A_n/(A_n\der_1+\cdots +A_n\der_n)$ has Gelfand-Kirillov dimension
 $n$. By Theorem \ref{FFI}  applied to the $A_n$-module $P_n$, we have
$$2n=\GK(A_n)\leq n ({\rm d}(A)+\max \{{\rm d}(A), 1\}),$$
hence ${\rm d }(A_n)\geq 1$, and so ${\rm d }(A_n)= 1$. $\Box $

\bigskip
{\sc Proof of the Bernstein's inequality (Theorem \ref{Berin})}.
\bigskip

Since $\GK (A_n)=2n$ and $\fdim (A_n)= 1$, Theorem \ref{FFI} gives
$\GK (M)\geq \frac{2n}{2}=n$.
 $\Box $

One also gets a short proof of the following result of Rentschler
and Gabriel.

\begin{theorem}\label{KdimAnRGab}
(\cite{Ren-Gab}). If char$\, K=0$ then the Krull dimension of the
Weyl algebra $A_n$ is
$$\Kdim (A_n)=n.
$$
\end{theorem}

{\it Proof}. Putting $\GK (A_n)=2n$ and $\fdim (A_n)= 1$ into the
second  formula  of Theorem \ref{SFI}
 we have $\Kdim (A_n)\leq \frac{2n}{2}=n$. The polynomial algebra $P_n=K[X_1, \ldots ,
X_n]$ is the subalgebra of $A_n$ such that $A_n$ is a free right
$P_n$-module.
   The map $I\ra A_n I$ from
 the set of  left ideals of the polynomial algebra  $P_n$ to the set of
 left ideals  of the Weyl algebra $A_n$ is injective, thus
 $n={\rm K.dim}(P_n)\leq \Kdim (A_n)$, and so $\Kdim (A_n)=n$. $\Box $


\section{Filter dimension of the ring of differential operators on a smooth
 irreducible affine algebraic variety (proof of Theorem \ref{fdif=1})
}\label{5fdiffhaz}

 Let $K$ be a field of characteristic $0$ and let the algebra $B$ be as in
 the previous section, i.e.
 $B$ is a finitely generated regular commutative algebra  which is a
 domain. We keep the notations of the previous section. Recall that
 the derivation ring $\Delta =\Delta (B) $ coincides with the
ring of differential operators ${\cal D} (B)$ (\cite{MR}, 15.5.6)
and is a simple finitely generated finitely partitive $K$-algebra
(\cite{MR}, 15.3.8, 15.1.21).
 We refer the reader to \cite{MR}, Chapter 15, for basic definitions.
 We aim to prove Theorem \ref{fdif=1}.

  Let $\{ B_i\}$ and $\{\D_i\}$ be standard finite dimensional filtrations
on $B$ and $\D $ respectively such that $B_i\subseteq \D_i$ for
all $i\geq 0$.
 Then the enveloping algebra $\D^e :=\D \t \D^{0}$ can be equipped with
 the  standard finite dimensional filtration $\{\D^e_i\}$ which is the
tensor product of the filtrations $\{ \D_i\}$ and $\{ \D_i^{0} \}$
of the algebras $\D $ and $\D^0$ respectively.

Then $B\simeq \D /\D {\rm Der}_KB$ is a simple left $\D $-module
(\cite{MR}, 15.3.8] with  $\GK (\D )$ $=2\GK (B)$ (\cite{MR},
15.3.2). By Theorem \ref{FFI},
$$ \fdim (\D )+\max \{ \fdim ( \D ), 1\}\geq \frac{\GK (\D )}{\GK (B)}=\frac{2\GK (B)}{\GK
(B)}=2,$$ hence $\fdim ( \D ) \geq 1$. It remains to prove the
opposite inequality. For,  we recall some properties of $\D $ (see
\cite{MR}, Ch. 15, for details).

 Given $0\neq c\in B$, denote by $B_c$ the localization of the algebra $B$
 at the powers of the element $c$, then $\Delta (B_c)\simeq \D (B)_c$ and the map
 $\D(B)
\ra \D (B)_c$, $d\ra d/1$, is injective  (\cite{MR}, 5.1.25).
There is a finite subset $\{ c_1, \ldots ,c_t\}$ of $B$ such that
the algebra  $\prod^t_{i=1}\;\Delta (B_{c_i})$
 is left and right faithfully flat over its subalgebra $\D  $,
\begin{equation}\label{4.1haz}
 \sum^t_{i=1}\; Bc_i=B \;\;{\rm (see\;the \; proof\;of\;15.2.13,
\cite{MR})}.
\end{equation}

For each $c=c_i$, ${\rm Der}_K(B_c)$ is a free $B_c$-module with a
basis  $\derj =\frac{\der }{\der x_j}$, $j=1,\ldots ,n$ for some
$x_1, \ldots ,x_n\in B$
 (\cite{MR}, 15.2.13). Note that the choice of the  $x_j$'th depends on the
 choice of the $c_i$. Then
 $$ \D (B)_c\simeq \D (B_c)=B_c < \der_1, \ldots ,
 \der_n> \supseteq K< x_1, \ldots , x_n, \der_1, \ldots ,
 \der_n> $$

Fix $c=c_i$. We aim to prove the following statement:
\begin{itemize}
\item  {\it there exist natural numbers} $a$, $b$, $\a $, {\it
and} $\b $ {\it such that for any } $0\neq d \in \D_k$ {\it there
exists }
 $$ w\in \D^e_{ak+b}:\;\;\;wd=c^{\a k +\b }.
\leqno(*)
$$
\end{itemize}

Suppose that we are done. Then one can choose the numbers
 $a$, $b$, $\a $,  and $\b $ such
that (*) holds for all $i=1, \ldots , t$. It follows from
(\ref{4.1haz}) that
$$ \sum^t_{i=1}\; f_ic_i=1 \;\; {\rm for \; some }\; f_i\in A.$$
Choose $\nu \in \mathbb{N} $: all $f_ic_i \in \D_{\nu }$ and set
$N(k)=\a k+\b $, then
$$1=(\sum^t_{i=1} \,f_ic_i)^{tN(k)}=\sum^t_{i=1} \,g_ic_i^{N(k)}=
\sum^t_{i=1} \, g_iw_id=wd,$$ where the $w_i$ are  from (*), i.e.
$w_i \in \D^e_{ak+b}$, $w_id=c^{N(k)}_i$. So, $w=\sum^t_{i=1}
\,g_iw_i\in $ $ \D^e_{\nu tN(k)+ak +b}$ and so
 $\fdim (\D ) \leq 1$, as required.

Fix $c=c_i$. By (\cite{MR}, 15.1.24) ${\rm Der}_K(B_c)\simeq {\rm
Der}_K(B)_c$ and ${\rm Der}_KB$ can be seen as a finitely
generated $B$-submodule of ${\rm Der}_K(B_c)$ (\cite{MR}, 15.1.7).

The algebra $B$ contains the polynomial subalgebra $P=K[x_1,
\ldots , x_n]$. The polynomial algebra $P$ has the natural
filtration $P=\cup_{i\geq 0}P_i$ by the total degree of the
variables. Fix a natural number $l$ such that $P_1\subseteq B_l$,
then $P_i\subseteq B_{li}$ for all $i\geq 0$. We denote by
$Q=K(x_1, \ldots , x_n)$ the field of fractions of $P$. The field
of fractions, say $L$, of the algebra $B$ has the same
transcendence degree $n$ as the field of rational functions $Q$.
The algebra $B$ is a finitely generated algebra, hence $L$ is a
finite field extension of $Q$ of dimension, say $m$, over $Q$. Let
$e_1, \ldots , e_m\in B$ be a $Q$-basis for the vector space $L$
over $Q$. Note that $L=QB$. One can find a natural number $\beta
\geq 1$ and a nonzero polynomial $p\in P_\beta $ such that
$$ \{ B_1, e_je_k\, | \, j,k=1, \ldots , m\}\subseteq
 \sum_{\alpha =1}^mp^{-1} P_\beta e_\alpha . $$ Then $B_k\subseteq
\sum_{j=1}^mp^{-2k}P_{2\beta k}e_j$ and $B_ke_i\subseteq
\sum_{j=1}^mp^{-3k}P_{3\beta k}e_j$ for all $k\geq 1$ and $i=1,
\ldots , m$. Let $0\neq d\in B_k$. The $m\times m$ matrix of the
bijective $Q$-linear map $L\ra L$, $x\mapsto dx$, with respect to
the basis $e_1, \ldots , e_m$ has entries from the set
$p^{-3k}P_{3\beta k}$. So, its characteristic polynomial
$$ \chi_d(t)=t^m+\alpha_{m-1}t^{m-1}+\cdots +\alpha_0$$
has coefficients in $p^{-3mk}P_{3m\beta k}$, and $\alpha_0\neq 0$
as $x\mapsto dx$ is a bijection. Now, 
\begin{equation}\label{pahaz}
P_{6m\beta k}\ni
p^{3mk}\alpha_0=p^{3mk}(-d^{m-1}-\alpha_{m-1}d^{m-2}-\cdots
-\alpha_1)d\in B_{4m\beta k}P_{3m\beta k}\subseteq B_{m\beta k(4+3
l)}d.
\end{equation}
Let $\d_1, \ldots , \d_t$ be a set of generators for the left
$B$-module $\Der_K(B)$. Then
$$ \der_i\in \sum_{j=1}^tc^{-l_1}B_{l_1}\d_j\;\; {\rm for }\;\;
i=1, \ldots , n,$$ for some natural number $l_1$. Fix a natural
number $l_2$ such that $\d_j(B_1)\subseteq B_{l_2}$ and $\d_j
(c)\in B_{l_2}$ for $j=1, \ldots , t$.  Then
$$\der^\alpha (B_k)\subseteq c^{-2| \alpha | (l_1+1)}B_{k+| \alpha |
(l_1+l_2)}\;\; {\rm for\; all}\;\; \alpha \in \mathbb{N}^n, \;
k\geq 1, $$ where $\alpha =(\alpha_1, \ldots , \alpha_n)$,
$|\alpha | = \alpha_1+\cdots +\alpha_n$, $\der^\alpha =
\der_1^{\alpha_1} \cdots \der_n^{\alpha_n}$. It follows from
(\ref{pahaz}) that one can find $\alpha \in \mathbb{N}^n$ such
that $| \alpha | \leq 6m\beta k$ and
$$ 1\in K^* \der^\alpha (p^{3mk}\alpha_0)\subseteq \der^\alpha
(B_{m\beta k(4+3l)}d)\subseteq c^{-12m\beta
(l_1+1)k}\D^e_{2(m\beta k(4+3l)+6m\beta k (l_1+l_2))}d$$ where
$K^*=K\backslash \{ 0\}$.
 Now (*) follows. $\Box $

In fact we have proved the following corollary.

\begin{corollary}
There exist natural numbers $a$ and $b$ such that for any $0\neq
d\in \D
 {}_k$ there exists an element $w\in \D^e_{ak+b}$ satisfying $wd=1$.
\end{corollary}


\section{Multiplicity for the filter dimension,
holonomic modules over simple finitely generated
algebras}\label{multholhaz}

In this section, we introduce a concept of multiplicity for the
filter dimension and a concept of holonomic module for (some)
finitely generated algebras. We will prove that a holonomic module
has finite length (Theorem \ref{holfhaz}). The multiplicity for
the filter dimension is a key ingredient in the proof.

First we recall the definition of multiplicity in the commutative
situation and then for certain non-commutative algebras (somewhat
commutative algebras).

{\bf Multiplicity in the commutative situation}. Let $B$ be a
commutative finitely generated $K$-algebra with a standard finite
dimensional filtration $F=\{ B_i\}$, and let $M$ be a finitely
generated $B$-module with a finite dimensional generating
subspace, say $M_0$, and with the standard filtration $\{
M_i=B_iM_0\}$ attached to it. Then there exists a polynomial
$p(t)=lt^d+\cdots \in \mathbb{Q}[t]$  with rational coefficients
of degree $d=\GK (M)$ such that $$ \dim_K (M_i)=p(i)\;\; {\rm for
\; all}\;\; i\gg 0.$$ The polynomial $p(t)$ is called the {\bf
Hilbert polynomial} of the $B$-module $M$. The Hilbert polynomial
does depend on the filtration $\{M_i\}$ of the module $M$ but its
leading coefficient $l$ does not. The number $e(M)=d!l$ is called
the {\bf multiplicity} of the $B$-module $M$. It is a natural
number which does  depend on the filtration $F$ of the algebra
$B$.

In the case when $M=B$ is the homogeneous coordinate ring of a
projective algebraic variety $X\subseteq \mathbb{P}^m$ equipped
with the natural filtration that comes from the grading of the
graded algebra $B$, the multiplicity is the {\bf degree} of $X$,
the number of points in which $X$ meets a general plane of
complementary degree in $\mathbb{P}^m$ ($K$ is an algebraically
closed field).

{\bf Somewhat commutative algebras}.  A $K$-algebra $R$ is called
a {\bf somewhat commutative algebra}  if it has a finite
dimensional filtration $R=\cup_{i\geq 0}\, R_i$ such that the
associated graded algebra $\gr \, R:=\oplus_{i\geq 0}\,
R_i/R_{i-1}$  is a commutative finitely generated  $K$-algebra
where $R_{-1}=0$ and $R_0=K$.
 Then the algebra  $R$ is a Noetherian  finitely generated algebra since $\gr \, R$ is so.
A finitely generated module over a somewhat commutative algebra
has the
 Gelfand-Kirillov dimension which is a natural number. We refer the reader to
the books \cite{KL, MR} for the properties of somewhat commutative
algebras.

{\it Definition}. For a  somewhat commutative algebra $R$ we
define the {\bf holonomic number},
$$ h_R:=\min \{ \GK (M)\; | \; M\neq 0\;{\rm is \; a \; finitely \;
generated}\; R-{\rm module}\}.$$

{\it Definition}. A finitely generated $R$-module $M$ is called a
{\bf holonomic} module if $\GK (M)=h_R$. In other words, a nonzero
finitely generated $R$-module is holonomic iff it has the least
 Gelfand-Kirillov dimension. If $h_R=0$ then every holonomic
$R$-module is finite dimensional and vice versa.

{\it Examples}. $1$.  The holonomic number of the Weyl algebra
$A_n$ is $n$. The polynomial algebra $K[X_1, \ldots , X_n]\simeq
A_n/\sum_{i=1}^n\, A_n \der_i $ with the natural action of the
ring of differential operators $A_n=K[X_1, \ldots , X_n, {\der
\over \der X_1}, \ldots , {\der \over \der X_n}]$ is
 a simple holonomic $A_n$-module.

 $2$. Let $X$ be a smooth irreducible affine algebraic variety
of dimension $n$. The ring of differential operators $\CD (X)$ is
a simple somewhat commutative algebra of Gelfand-Kirillov
dimension $2n$ with  holonomic number $h_{\CD (X)}=n$.
 The algebra $\OO (X)$ of regular functions of the variety $X$ is a simple
$\CD (X)$-module with respect to the natural action of the algebra
$\CD (X)$. In more detail, $\OO (X)\simeq \CD (X)/\CD
(X)\Der_K(\OO (X))$ where $\Der_K(\OO (X))$ is the $\OO
(X)$-module of derivations of the algebra $\OO (X)$.

Let $R=\cup_{i\geq 0}\, R_i$  be a somewhat commutative algebra.
The associated graded
 algebra $\gr \, R$ is a commutative affine algebra. Let us choose homogeneous algebra
 generators of the algebra $\gr \, R$, say $y_1, \ldots , y_s$, of graded degrees
$1\leq k_1, \ldots , k_s$ respectively (that is $y_i \in
R_{k_i}/R_{k_i-1}$).
 A filtration $\G =\{ \G_i , i\geq 0\}$ of an $R $-module
$M=\cup^{\infty }_{i=0}\, \G_i$ is called {\bf good} if the
associated graded
 $\gr \, R $-module $\gr_\G M:=\oplus_{i\geq 0}\,\G_i/ \G_{i-1}$ is finitely generated.
 An $R $-module $M$ has
 a good filtration iff it is finitely generated, and if $\{ \G_i \}$ and
$\{ \Omega_i \}$ are two good filtrations of $M$, then there
exists a natural number $t$
 such that $\G_i \subseteq \Omega_{i+t}$ and $\Omega_i \subseteq \G_{i+t}$ for all $i$.
 If an $R $-module $M$ is finitely generated and $M_0$ is a
finite dimensional generating subspace of $M$, then the standard
filtration
 $\{ \G_i =R_iM_0\}$ is good (see \cite{Bj, Bor, Co, KL, MR, NVO}  for
details).  The following two Lemmas are well-known by specialists
(see their proofs, for example, in \cite{Bav 4}, Theorem 3.2 and
Proposition 3.3 respectively).

\begin{lemma}\label{kpolynomials}
Let  $ R =\cup_{i\geq 0}\;R_i $ be a somewhat commutative algebra,
$k={\rm lcm}(k_1, \ldots , k_s)$, and let $M$ be a finitely
generated $R $-module with good filtration $\G =\{ \G_i \}$.
\begin{enumerate}
\item There exist $k$ polynomials $\g_0 , \ldots ,\g_{k-1}\in
\mathbb{Q} [t]$
 with  coefficients from $[k^{\GK(M)} \GK(M)!]^{-1} \mathbb{Z} $ such that

$$ \dim\,  \G_i=\g_j(i)\;\;{\rm for \;all}\; i\gg 0 \; {\rm and}\;
j\equiv i \, ({\rm mod}\; k).$$

\item The polynomials $\g_j$ have the same degree $\GK (M)$ and
the same leading coefficient $e(M)/\GK (M)!$ where $e(M)$ is
called the {\bf multiplicity}
 of $M$. The multiplicity $e(M)$ does not depend on the choice of the good
 filtration $\G $. $\Box $
\end{enumerate}
\end{lemma}

{\it Remark}. A finitely generated $R$-module $M$ has $e(M)=0$ iff
$\dim_K (M)< \infty $.

\begin{lemma}\label{multiplicityexact}
Let $0\ra N\ra M\ra L\ra 0$ be an exact sequence of modules over a
somewhat commutative
 algebra $R$. Then $\GK (M)=\max \{ \GK (N), \GK (L) \}$, and if
 $\GK (N)= \GK (M)=\GK (L)$ then $e(M)=e(N)+e(L)$. $\Box $
\end{lemma}

\begin{corollary}\label{holmodfinlength} 
Let the algebra $R$ be as in Lemma \ref{kpolynomials} with
holonomic number $h>0$.
\begin{enumerate}
\item Let $M$ be a holonomic $R$-module with multiplicity $e(M)$.
The $R$-module $M$ has finite length $\leq e(M) k^h$. \item Every
nonzero submodule or factor module of a holonomic $R$-module is a
holonomic
 module.
\end{enumerate}
\end{corollary}

{\it Proof}. This follows directly from Lemma
\ref{multiplicityexact}. $\Box $

{\bf Multiplicity}.  Let $f$ be a function from $\mathbb{N}$ to
$\mathbb{R}_+=\{ r\in \mathbb{R}: \, r\geq 0\}$, the {\bf leading
coefficient} of
 $f$ is a non-zero limit (if it exists)
$$\lc (f)=\lim \frac{f(i)}{i^d}\not= 0, \,\,\,i\ra \infty ,$$
where $d=\g (f)$. If $d\in \mathbb{N} $, we define the {\bf
multiplicity} $e(f)$ of $f$ by
$$e(f)=d!\,\lc (f).$$
The factor $d! $ ensures that the multiplicity $e(f)$ is a
positive integer in some important cases. If
$f(t)=a_dt^d+a_{d-1}t^{d-1}+\cdots + a_0$ is a polynomial of
degree $d$ with real coefficients then $\lc (f)=a_d$ and
$e(f)=d!a_d$.

\begin{lemma}\label{lchaz}
 Let $A$ be a finitely generated algebra equipped with a standard
finite dimensional filtration $F=\{ A_i\}$ and $M$ be a finitely
generated $A$-module with generating finite dimensional subspaces
$M_0$ and $N_0$.

1. If  $\lc (\nFM0 )$ exists then so does $\lc (\nFN0 )$, and $\lc
(\nFM0 )=\lc (\nFN0 )$.

2. If   $\lc (\dim \, A_iM_0 )$ exists then so does $\lc (\dim \,
A_iN_0)$, and $\lc (\dim \, A_iM_0)=\lc (\dim \, A_iN_0)$.
\end{lemma}

{\it Proof}. 1. The module $M$ has two filtrations $\{M_i=A_iM_0
\}$ and $\{N_i=A_iN_0 \} $.  Let $\nu  =\nFM0 $    and $\mu =\nu
{}_{F, N_0} $. Choose a natural number $s$ such that $M_0\subseteq
N_s$ and $N_0\subseteq M_s$, so $N_i\subseteq M_{i+s}$ and
$M_i\subseteq N_{i+s}$ for
 $i\geq 0$. Since $M_0\subseteq A_{\nu (i+s)}N_{i,gen}$ for each $i$  and $N_0
\subseteq A_sM_0,$ we have $N_0\subseteq A_{\nu
(i+s)+s}N_{i,gen},$ hence, $\mu (i)\leq \nu (i+s)+s$. By symmetry,
$\nu (i)\leq \mu (i+s)+s$, so if  $\lc (\mu )$ exists then so does
$\lc (\nu )$ and  $\lc (\mu )=\lc (\nu )$.

2. Since $\dim \,  N_i \leq \dim \, M_{i+s}$ and $\dim \, M_i\leq
\dim \, N_{i+s}$ for $i\geq 0$, the statement is clear. $\Box $

 Lemma \ref{lchaz} shows that the leading coefficients of the functions
$\dim \, A_iM_0$ and $\nFM0 $ (if exist) do not depend on the
choice of the generating subspace $M_0$. So, denote them by
$$l(M)=l_F(M)\;\; {\rm and}\;\; L(M)=L_F(M)$$
 respectively (if they exist). If
$\GK (M)$ (resp. $\fdim (A)$) is a natural number, then we denote
by $e(M)=e_F(M)$ (resp. $E(M)=E_F(M)$) the multiplicity of the
 function $\dim \, A_iM_0$ (resp. $\nFM0 $).

We denote by $L(A)=L_F(A)$ the leading coefficient $L_F({}_{A\t
A^0}A)$ of the return function ${\nu }_{F\t F^0,K}$ of the $A\t
A^0$-module $A$.

{\bf Holonomic modules}.  {\it Definition}. Let $A$ be a finitely
generated $K$-algebra, and $h_A$ be its holonomic number. A
nonzero finitely generated $A$-module $M$ is called a {\bf
holonomic} $A$-module if $\GK (M)=h_A$. We denote  by $\hol (A)$
the set of all the holonomic $A$-modules.

Since the holonomic number is an infimum it is not clear at the
outset that there will be modules which achieve this dimension.
Clearly, $\hol(A)\neq \emptyset $ if the Gelfand-Kirillov
dimension of every finitely generated $A$-module is a natural
number.

A nonzero submodule or a factor module of a holonomic is a
holonomic module (since the Gelfand-Kirillov dimension of a
submodule or a factor module does not exceed the Gelfand-Kirillov
of the module). If, in addition, the finitely generated algebra
$A$ is left Noetherian and finitely partitive then each holonomic
$A$-module $M$ has finite length and each simple sub-factor of $M$
is a holonomic module.

 Let us consider algebras $A$ having the following properties:

\begin{itemize}

\item (S) $A$ {\em is a simple finitely generated infinite
dimensional algebra.}

\item (N) {\em There exists a standard finite dimensional
filtration}
 $F=\{A_i\}$ {\em of the algebra} $A$ {\em such that  the
associated graded algebra} gr $A:=$ $\oplus {}_{i\geq
0}\,A_i/A_{i-1}$, $A_{-1}=0$,
 {\em is left Noetherian.}

\item (D) $\GK (A)<\infty $, $\fd  (A)<\infty $, {\em both}
$l(A)=l_F(A)$
 {\em and} $L(A)=
L_F(A)$ {\em exist.}

\item (H) {\em For every holonomic} $A$-{\em module} $M$ {\em
there exists}
 $l(M)=l_F(M)$.

\end{itemize}

In many cases we use the weaker form of the condition (D).
\begin{itemize}

\item (D') $\GK (A)<\infty $, $d=\fd (A)<\infty $, {\em there
exist}
 $l(A)=l_F(A)$ {\em and
a positive number} $c>0$ {\em such that} $\nu (i)\leq ci^d$ {\em
for}
 $i\gg 0$ {\em where}
$\nu $ {\em is the return function} $\nu {}_{F\t F^0,K}$ {\em of
the left } $A\t A^0$-{\em module}  $A$.

\end{itemize}

It follows from (N) that $A$ is a  left Noetherian algebra.

\begin{lemma}
(\cite{Bavcafd})
\begin{enumerate}
\item The Weyl algebra $A_n$ over a field of  characteristic zero
with the standard finite dimensional filtration $F=\{A_{n,i}\}$
associated with the canonical generators satisfies the conditions
(S), (N), (D), (H). The return function $\nu_F(i)=i$ for $i\geq
0$, and so the leading coefficient of $\nu_F$ is $L_F(A_n)=1$.
\item $\nu_{G,K}(i)=i$ for $i\geq 0$ and $L_G(P_n)=1$ where
$\nu_{G,K}$ is the return function of  the $A_n$-module
$P_n=K[X_1, \ldots ,X_n]=A_n/(A_n{\der}_1+\cdots +A_n{\der }_n)$
with the usual filtration $G=\{P_{n,i}\}$ of the polynomial
algebra.
\end{enumerate}
\end{lemma}

{\it Proof}. $1$. The only fact that we need to prove is that
$\nu_F(i)=i$ for $i\geq 0$. We keep the notation of Theorem
\ref{fdWeyl=1}. In the proof of Theorem \ref{fdWeyl=1} we have
seen that $\nu_F(i)\leq i$ for $i\geq 0$. It remains to prove the
reverse inequality.

Each element $u$ in $A_n$ can be written in a unique way as a
finite sum $u=\sum \,\lab\Xa\db$ where  $\lab \in K$ and $\Xa $
denotes the  monomial $X^{\a {}_1}_1 \cdots X^{{\a }_n}_n$ and
similarly $\db $ denotes the  monomial $\partial {}^{\beta
{}_1}_1\cdots \partial {}^{\beta {}_n}_n$. The element $u$ belongs
to $A_{n,m}$ iff $|\a |+|\b |\leq m$, where $|\a |=\a {}_1+ \cdots
+ \a {}_n$. If $\a \in K[X_1,\ldots ,X_n]$, then
$${\partial }^m_i \alpha  ={\sum }^m_{j=0}\,{m\choose j}
\frac{\partial^j\alpha }{\partial X^j_i} {\partial
}^{m-j}_i,\,\,\,m\in \mathbb{N} .$$ It follows that for any $v\in
\sum \,A_{n,i}\t A_{n,j}^0$,
 $i+j<m$, the element $vX^m_1=\sum \,\lab \Xa \db $ has the coefficient
${\lambda }_{0, 0}=0$, hence it could not be a non-zero scalar,
and so  $\nu (i)\geq i$ for all $i\geq 0$. Hence $\nu (i)=i$ all
$i\geq 0$ and then $L_F(A_n)=1$.

$2$. The standard filtration of the $A_n$-module $P_n$ associated
with the generating subspace $K$ coincides with the usual
filtration of the polynomial algebra $P_n$. Since
$\der_j(P_{n,i})\subseteq P_{n,i-1}$ for all $i\geq 0$ and $j$,
$\nu_{G,K}(i)\leq i$ for $i\geq 0$. Using the same arguments as
above we see that for any $u\in \sum_{j=0}^{i-1}A_{n,j}\t
A^0_{n,i-j-1}$ the element $uX^i_1$ belongs to the ideal of $P_n$
generated by $X_1$, hence, $\nu_{G,K} (i)\geq i$, and so
$\nu_{G,K} (i)= i$ for all $i\geq 0$ and $L_G(P_n)=1$. $\Box $

\begin{theorem}\label{63haz}
 (\cite{Bavcafd}) Assume that an algebra $A$ satisfies the conditions (S), (H), (D) resp. (D')
for some standard finite dimensional filtration $F=\{A_i\}$ of
$A$. Then for every holonomic $A$-module $M$ its leading
coefficient is bounded from below by a nonzero constant:
$$ l(M)\geq \sqrt{\frac{l(A)}{(L(A)L'(A))^{h_A}} },$$ where

$$L'(A)=\left\{%
\begin{array}{ll}
    L(A), & \hbox{if $\fdim (A)>1$,} \\
   L(A)+1 , & \hbox{if $\fdim (A)=1$,} \\
    1, & \hbox{if $\fdim (A)<1$.} \\
\end{array}%
\right.$$
 resp.
$$ l(M)\geq \sqrt{\frac{l(A)}{(c(c+1))^{h_A}} }.$$
\end{theorem}

{\it Proof}. Let $M_0$ be a generating finite dimensional subspace
of $M$ and $\{M_i=A_iM_0\}$ be the  standard finite dimensional
filtration of $M$. In the
 proof of Theorem \ref{FFI} we proved that
$\dim \, A_i\leq \dim \, M_{\l (i)}\dim \, M_{\l (i)+i}$ for
 $i\geq 0$ where $\l $ is the left return function of the algebra  $A$ associated
with the filtration $F$. Since $\l (i)\leq \nu (i)$ for $i\geq 0$
we have  $\dim \, A_i\leq \dim \, M_{\nu (i)}\dim \, M_{\nu (i)
+i}$, hence, if (D) holds then
$$l(A)i^{\GK (A)}+\cdots \leq l^2(M)(L(A)L'(A))^{\GK (M)}
i^{\GK (M)(\fd (A)+\max \{\fd (A), 1\})}+\cdots ,$$ where three
dots denote  smaller terms.

If (D') holds then
$$l(A)i^{\GK (A)}+\cdots \leq l^2(M)(c(c+1))^{\GK (M)}
i^{\GK (M)(\fd (A)+\max \{\fd (A), 1\})}+\cdots .$$

The module $M$ is holonomic, i.e. $\GK (A)= \GK (M)(\fd (A)+\max
\{\fd (A),1\})$. Now, comparing the ``leading" coefficients in the
inequalities above we finish the proof. $\Box $

Let $A$ be as in Theorem \ref{63haz}. We attach to the algebra $A$
two positive numbers $c_A$ and $c'_A$ in the cases (D) and (D')
respectively:

$$c_A=\sqrt{\frac{l(A)}{(L(A)L'(A))^{h_A}} } \;\; {\rm and }\;\;
(c'_A)=\sqrt{\frac{l(A)}{(c(c+1))^{h_A}} }.$$

\begin{corollary}\label{67haz}
Assume that an algebra $A$ satisfies the conditions (S), (N), (H),
(D) or (D'). Let $0\ra N\ra M \ra L\ra 0$ be an exact sequence of
nonzero finitely generated $A$-modules. Then $M$ is holonomic if
and only if  $N$ and $L$ are holonomic, in that case
$l(M)=l(N)+l(L)$.
\end{corollary}

{\it Proof }. The algebra $A$ is left Noetherian, so the module
$M$ is finitely generated iff both $N$ and $L$ are so. The proof
of Proposition 3.11 (\cite{MR}, p. 295) shows that we can choose
finite dimensional generating subspaces
 $N_0$, $M_0$, $L_0$ of the modules $N$, $M$, $L$ respectively such that the
sequences
$$0\ra N_i=A_iN_0\ra M_i=A_iM_0\ra L_i=A_iL_0\ra 0$$
 are exact for all $i$, hence, $\dim\,  M_i=\dim \, N_i +\dim \, L_i$ and the results
follow. $\Box $

\begin{theorem}\label{holfhaz}
(\cite{Bavcafd}) Suppose that the conditions  (S), (N), (H), (D)
(resp. (D')) hold. Then each holonomic $A$-module $M$ has finite
length which is less or equal to $l(M)/c_A$ (resp. $l(M)/c'_A$).
\end{theorem}

{\it Proof}. If $M=M_1\supset M_2\supset \ldots \supset M_m\supset
M_{m+1}=0$ is a chain of distinct submodules, then by Corollary
\ref{67haz} and Theorem \ref{63haz}
$$l(M)={\sum }^m_{i=1}\,l(M_i/M_{i+1})\geq mc_A,\,\,\,({\rm resp.}\,\,\,
l(M)\geq mc'_A),$$
 thus $m\leq l(M)/c_A$ (resp. $m\leq l(M)/c'_A$). $\Box $


\section{Filter dimension and commutative subalgebras of simple finitely generated
 algebras and their division
algebras}\label{hazmaxissPA}

In this section, using the first and the second filter
inequalities, we obtain $(i)$ an {\em upper} bound for the
Gelfand-Kirillov dimension of commutative subalgebras of simple
finitely generated infinite dimensional algebras (Theorem
\ref{GKcsuba}), and $(ii)$ an {\em upper} bound for the
transcendence degree of  subfields of quotient division rings of
(certain) simple finitely generated infinite dimensional algebras
(Theorems \ref{trdsubf} and \ref{c1trds}).

For certain classes of algebras and their division algebras the
maximum Gelfand-Kirillov dimension/transcendence degree over the
commutative subalgebras/subfields were found in
\cite{AmitsurPAMS58}, \cite{GK66}, \cite{Mak-Limcom},
\cite{Joseph72HP}, \cite{JosLN74}, \cite{JosgenQ},
\cite{Amitsur-Small78}, and \cite{Resco79}.

Recall that
\begin{eqnarray*}
{\rm the \; Gelfand-Kirillov \; dimension}\; \GK (C)& =& {\rm the
\; Krull \; dimension}\; \Kdim (C)\\
 &=&   {\rm the \;  transcendence \; degree} \;
\trdeg_K (C)
\end{eqnarray*}
for every commutative finitely generated algebra $C$ which is a
domain.

{\bf An upper bound for the Gelfand-Kirillov dimensions of
commutative subalgebras of simple finitely generated algebras}.

 \begin{proposition}\label{drGK}
Let $A$ and $C$ be  finitely generated algebras such that $C$ is a
commutative domain with field of fractions $Q$, $B:=C\t A$, and
$\CB :=Q\t A$. Let $M$ be a finitely generated $B$-module such
that  $\CM :=\CB \t_BM\neq 0$. Then $ \GK ({}_BM)\geq \GK_Q({}_\CB
\CM )+\GK (C)$.
\end{proposition}

{\it Remark}. $\GK_Q$ stands for the Gelfand-Kirillov dimension
over the field $Q$.

{\it Proof}. Let us fix  standard filtrations  $\{ A_i\}$ and $\{
C_i\}$ for the algebras $A$ and $C$ respectively. Let $h(t)\in
\mathbb{Q}[t]$ be the {\em Hilbert polynomial} for the algebra
$C$, i.e. $\dim_K(C_i)=h(i)$ for $i\gg 0$. Recall that $\GK
(C)=\deg_t(h(t))$. The algebra $B$ has a standard filtration $\{
B_i\}$ which is the tensor product of the standard filtrations
$\{C_i\}$ and $\{ A_i\}$  of the algebras $C$ and $A$, i.e.
$B_i:=\sum_{j=0}^iC_j\t A_{i-j}$. By the assumption, the
$B$-module $M$ is finitely generated, so $M=BM_0$ where $M_0$ is a
finite dimensional generating subspace for $M$. Then the
$B$-module $M$ has a standard filtration $\{ M_i:=B_iM_0\}$. The
$Q$-algebra $\CB $ has a standard (finite dimensional over $Q$)
 filtration  $\{ \CB_i:=Q\t A_i\}$, and the $\CB $-module $\CM $
has a standard (finite dimensional over $Q$) filtration $\{
\CM_i:=\CB_iM_0'=QA_iM_0'\}$ where $M_0'$ is the image of the
vector space $M_0$ under the $B$-module homomorphism $M\ra \CM$,
$m\mapsto m':=1\t_Bm$.

For each $i\geq 0$, one can  fix a $K$-subspace, say $L_i$, of
$A_iM_0'$ such that $\dim_Q(QA_iM_0')=\dim_K(L_i)$. Now,
$B_{2i}\supseteq C_i\t A_i$ implies $\dim_K(B_{2i}M_0)\geq
\dim_K((C_i\t A_i)M_0)$, and $((C_i\t A_i)M_0)'\supseteq C_iL_i$
implies $\dim_K(((C_i\t A_i)M_0)')\geq
\dim_K(C_iL_i)=\dim_K(C_i)\dim_K(L_i)=\dim_K(C_i)\dim_Q(\CM_i)$.
It follows that
\begin{eqnarray*}
\GK ({}_BM)&=&\g (\dim_K(M_i))\geq \g (\dim_K(M_{2i}))=\g
(\dim_K(B_{2i}M_0))\geq \g (\dim_K((C_i\t A_i)M_0))\\
&\geq & \g (\dim_K(((C_i\t A_i)M_0)') \geq \g
(\dim_K(C_i)\dim_Q(\CM_i))\\
&=&\g (\dim_K(C_i))+ \g (\dim_Q(\CM_i))\;\; ({\rm since  }\;\; \g
(\dim_K(C_i))=h(i), \;\; {\rm for }\;\; i\gg 0)\\
&=& \GK (C)+\GK_Q({}_\CB \CM ).\;\; \Box
\end{eqnarray*}

Recall that $d=\fd, \lfd$.  A $K$-algebra $A$  is called {\em
central} if its centre $Z(A)=K$.
\begin{theorem}\label{GKcsuba}
(\cite{dimcom}) Let $A$ be a central simple finitely generated
$K$-algebra of Gelfand-Kirillov dimension $0<n<\infty $ (over
$K$). Let $C$ be a commutative subalgebra of $A$. Then
$$ \GK (C)\leq \GK (A)\left( 1-\frac{1}{f_A+\max \{
f_A, 1\}} \right)$$ where $f_A:= \max \{ \fdim_{Q_m}(Q_m\t A)\, |
\, 0\leq m\leq  n\}$, $Q_0:=K$, and $Q_m:=K(x_1, \ldots , x_m)$ is
a rational function field in indeterminates $x_1, \ldots , x_m$.
\end{theorem}

{\it Proof.} Let $P_m=K[x_1, \ldots , x_m]$ be a polynomial
algebra over the field $K$. Then $Q_m$ is its field of fractions
and $\GK (P_m)=m$. Suppose that $P_m$ is a subalgebra of $A$. Then
$m=\GK (P_m)\leq \GK (A)=n$. For each $m\geq 0$, $Q_m\t A$ is a
central simple $Q_m$-algebra (\cite{MR}, 9.6.9) of
Gelfand-Kirillov dimension (over $Q_m$) $\GK_{Q_m}(Q_m\t A)=\GK
(A)>0$, hence
\begin{eqnarray*}
\GK (A)&=& \GK({}_AA_A)\geq \GK({}_AA_{P_m})=\GK({}_{ P_m\t A
}A)\;\;\;
(P_m \;\; {\rm is \; commutative})\\
&\geq & \GK_{Q_m}({}_{Q_m\t A}(Q_m\t_{P_m}A))+\GK (P_m) \;\;\;
({\rm
Lemma} \; \ref{drGK})\\
 &\geq & \frac{\GK (A)}{\fdim_{Q_m}(Q_m\t A)+\max \{ \fdim_{Q_m}(Q_m\t A),
 1\}}+m\;\;\; ({\rm Theorem }\; \ref{FFI}).
\end{eqnarray*}
Hence,
$$ m\leq \GK (A)\left( 1-\frac{1}{\fdim_{Q_m} (Q_m\t A)+\max \{
\fdim_{Q_m} (Q_m\t A), 1\}} \right)\leq \GK (A), $$ and so
$$ \GK (C)\leq \GK (A)\left( 1-\frac{1}{f_A+\max \{
f_A, 1\}} \right). \;\; \Box $$

 As a
consequence we have a short proof of the following well-known
result.

\begin{corollary}\label{CDXn}
Let $K$ be an algebraically closed  field of characteristic zero,
$X$ be a smooth irreducible affine  algebraic variety of dimension
$n:=\dim (X)>0$, and $C$ be a commutative subalgebra of the ring
of differential operators $\CD (X)$. Then $\GK (C)\leq n$.
\end{corollary}

{\it Proof}. The algebra $\CD (X)$ is central since $K$ is an
algebraically closed field of characteristic zero  \cite{MR}, Ch.
15.  By Theorem \ref{fdif=1}, $f_{\CD (X)}=1$, and then, by
Theorem \ref{GKcsuba},
$$ \GK (C)\leq 2n (1- \frac{1}{1+1})=n. \; \Box $$

{\it Remark}. For the ring of differential operators $\CD (X)$ the
upper bound in Theorem \ref{GKcsuba} for the Gelfand-Kirillov
dimension  of commutative subalgebras of $\CD (X)$ is an {\em
exact}  upper  bound since as we mentioned above the  algebra $\OO
(X)$ of regular functions on $X$ is a commutative subalgebra of
$\CD (X)$ of Gelfand-Kirillov dimension $n$.

{\bf An upper bound for the transcendence degree  of
 subfields of quotient division  algebras of simple finitely generated algebras}.

 Recall that the transcendence degree $\trdeg_K(L)$ of a field
extension $L$ of a field $K$  coincides with the Gelfand-Kirillov
dimension $\GK_K(L)$, and, by the {\bf Goldie's Theorem},  a left
Noetherian algebra $A$ has a quotient  algebra  $D=D_A$ (i.e.
$D=\S1 A$ where $S$ is the set of {\em regular} elements $=$ the
set of non-zerodivisors of $A$). As a rule, the quotient  algebra
$D$ has {\em infinite} Gelfand-Kirillov dimension and is {\em not}
a finitely generated algebra ({\em eg}, the quotient division
algebra $D(X)$ of the ring of differential operators $\CD (X)$ on
{\em each} smooth irreducible affine  algebraic variety $X$ of
dimension $n>0$ over a field $K$ of characteristic zero contains a
{\em noncommutative free} subalgebra since $D(X)\supseteq
D(\mathbb{A}^1)$ and the {\em first Weyl division algebra}
$D(\mathbb{A}^1)$ has this property \cite{Mak-LimFree}). So, if we
want to find an upper bound for the transcendence degree of
subfields in the quotient algebra $D$ we can not apply Theorem
\ref{GKcsuba}. Nevertheless, imposing some natural (mild)
restrictions on the algebra $A$ one can obtain exactly the same
upper bound for the transcendence degree of subfields in the
quotient algebra $D_A$ as the upper bound for the Gelfand-Kirillov
dimension of commutative subalgebras in $A$.

\begin{theorem}\label{trdsubf}
(\cite{dimcom}) Let $A$ be a simple finitely generated $K$-algebra
such that $0<n:=\GK (A)<\infty $, all the algebras $Q_m\t A$,
$m\geq 0$, are
 simple finitely partitive algebras where $Q_0:=K$, $Q_m:=K(x_1,
\ldots , x_m)$ is a rational function field and, for  each $m\geq
0$, the Gelfand-Kirillov dimension (over $Q_m$)  of every finitely
generated $Q_m\t A$-module is a natural number. Let $B=\S1 A$ be
the  localization of the algebra $A$ at a left Ore subset $S$ of
$A$. Let $L$ be a (commutative) subfield of the algebra $B$ that
contains $K$. Then
$$ \trdeg_K(L)\leq \GK (A)\left( 1-\frac{1}{f_A+\max \{
f_A, 1\}} \right)$$ where $f_A:=\max \{ \fdim_{Q_m}(Q_m\t A)\, |
\, 0\leq m \leq  n\}$.
\end{theorem}

{\it Proof}. It follows immediately from  a definition  of the
Gelfand-Kirillov dimension that $\GK_{K'}(K'\t C)=\GK (C)$ for any
$K$-algebra $C$ and any field extension $K'$ of $K$. In
particular, $\GK_{Q_m}(Q_m\t A)=\GK (A)$ for all $m\geq 0$. By
Theorem \ref{SFI},
$$ \Kdim (Q_m\t A)\leq \GK (A)\left( 1-\frac{1}{\fdim_{Q_m} (Q_m\t A)+\max \{
\fdim_{Q_m} (Q_m\t A), 1\}} \right).$$ Let $L$ be a subfield of
the algebra $B$ that contains $K$. Suppose that $L$ contains a
rational function field (isomorphic to) $Q_m$ for some $m\geq 0$.
\begin{eqnarray*}m&=& \trdeg_K(Q_m)\leq \Kdim (Q_m\t Q_m)\\
&\leq & \Kdim (Q_m\t B) \; ({\rm by}\; \cite{MR}, \, 6.5.3\;\;
{\rm since}\;
Q_m\t B \; {\rm is \; a\; free}\;\;  Q_m\t Q_m-{\rm module}) \\
&=& \Kdim (Q_m\t \S1 A)=\Kdim (\S1 (Q_m\t A))\\
&\leq & \Kdim (Q_m\t A)
\; \; ({\rm by}\; \cite{MR}, \, 6.5.3.(ii).(b))\\
&\leq & \GK (A)\left( 1-\frac{1}{\fdim_{Q_m} (Q_m\t A)+\max \{
\fdim_{Q_m} (Q_m\t A), 1\}} \right)\leq \GK (A).
\end{eqnarray*}
Hence
$$\trdeg_K(L)\leq \GK (A)\left( 1-\frac{1}{f_A+\max \{
f_A, 1\}} \right). \;\; \Box $$

Recall that every
 somewhat commutative algebra $A$ is a Noetherian  finitely generated finitely
  partitive algebra of finite Gelfand-Kirillov dimension, the Gelfand-Kirillov
  dimension of every finitely generated $A$-modules is an integer, and
 ({\bf Quillen's Lemma}): {\em the ring ${\rm End}_A(M)$ is algebraic
 over} $K$ (see
  \cite{MR}, Ch. 8 or \cite{KL} for details).

\begin{theorem}\label{c1trds}
(\cite{dimcom}) Let $A$ be a central simple  somewhat commutative
infinite dimensional $K$-algebra  and let $D=D_A$ be its quotient
algebra. Let $L$ be a subfield of $D$ that contains $K$. Then the
transcendence degree of the field $L$ (over $K$)
$$\trdeg_K (L)\leq \GK (A)\left( 1-\frac{1}{f_A+\max \{
f_A, 1\}} \right)$$ where $f_A:= \max \{ \fdim_{Q_m}(Q_m\t A)\, |
\, 0\leq m\leq \GK(A)\}$.
\end{theorem}

{\it  Proof.} The algebra $A$ is a somewhat commutative algebra,
so it has a finite dimensional filtration $A=\cup_{i\geq 0}A_i$
such that the associated graded algebra is a commutative finitely
generated algebra. For each integer $m\geq 0$, the $Q_m$-algebra
$Q_m\t A=\cup_{i\geq 0}Q_m\t A_i$ has the finite dimensional
filtration
 (over $Q_m$) such that the associated graded algebra $\gr (Q_m\t
 A)=\oplus_{i\geq 0}Q_m\t A_i/Q_m\t A_{i-1}\simeq Q_m\t \gr (A)$
 is a commutative finitely generated $Q_m$-algebra. So, $Q_m\t A$
 is a somewhat commutative $Q_m$-algebra.

 By the assumption $\dim_K(A)=\infty$, hence $\dim_K(\gr
 (A))=\infty $ which implies $\GK (\gr (A))>0$, and so $\GK (A)>0$
 (since $\GK (A)=\GK (\gr (A))$). The algebra $A$ is a central
 simple $K$-algebra, so $Q_m\t A$ is a central simple
 $Q_m$-algebra (\cite{MR}, 9.6.9). Now,  Theorem
 \ref{c1trds} follows from Theorem \ref{trdsubf} applied to $B=D$.
 $\Box $

\begin{theorem}\label{c2trds}
Let $K$ be an algebraically closed field of characteristic zero,
$\CD (X)$ be the ring of differential operators on a smooth
irreducible affine  algebraic variety $X$ of dimension $n>0$, and
$D(X)$ be the quotient division ring for $\CD (X)$. Let $L$ be a
(commutative) subfield of $D(X)$ that contains $K$. Then
$\trdeg_K(L)\leq n$.
\end{theorem}

{\it Remark}. This inequality is, in fact, an {\em exact} upper
bound for the transcendence degree of subfields in $D(X)$ since
the field of fractions $Q(X)$ for the algebra $\OO (X)$ is a
commutative subfield  of the division ring $D(X)$ with
$\trdeg_K(Q(X))=n$.

{\it Proof.} Since $Q_m\t \CD_K (\OO (X))\simeq \CD_{Q_m} (Q_m\t
\OO (X))$ and $\fdim (\CD (Q_m\t \OO (X)))=1$ for all $m\geq 0$ we
have $f_{\CD (X)}=1$. Now, Theorem \ref{c2trds} follows from
Theorem \ref{c1trds}, $$ \trdeg_K(L)\leq 2n(1-\frac{1}{1+1})=n. \;
\; \Box
$$

Following \cite{JosgenQ} for a $K$-algebra $A$ define the {\bf
commutative dimension}
$$ \Cdim (A):=\max \{ \GK (C)\, | \; C \;\; {\rm is \; a \;
commutative\; subalgebra \; of }\; A\}.$$ The commutative
dimension $\Cdim (A)$ (if finite) is the largest non-negative
integer $m$ such that the algebra $A$ contains a polynomial
algebra in $m$ variables (\cite{JosgenQ}, 1.1, or \cite{MR},
8.2.14). So, $\Cdim (A) =\mathbb{N}\cup \{ \infty \}$. If $A$ is a
subalgebra of $B$ then $\Cdim (A)\leq \Cdim (B)$.

\begin{corollary}\label{DXYnmnot}
Let $X$ and $Y$ be smooth irreducible affine  algebraic varieties
of dimensions $n$ and $m$ respectively, let $D(X)$ and $D(Y)$ be
quotient division rings for the rings of differential operators
$\CD (X)$ and $\CD (Y)$. Then there is no $K$-algebra embedding
$D(X)\ra D(Y)$ if $n>m$.
\end{corollary}

{\it Proof}. By Theorem \ref{c2trds}, $\Cdim (D(X))=n$ and $\Cdim
(D(Y))=m$. Suppose that there is a $K$-algebra embedding $D(X)\ra
D(Y)$. Then $n=\Cdim (D(X))\leq \Cdim (D(Y))=m$.
 $\Box $

For the Weyl algebras $A_n=\CD ( \mathbb{A}^n)$ and  $A_m=\CD (
\mathbb{A}^m)$ the result above was proved by Gelfand and Kirillov
in \cite{GK66}. They introduced a new  invariant of an  algebra
$A$, so-called the {\em (Gelfand-Kirillov) transcendence degree}
$\GKtrdeg (A)$, and proved that $\GKtrdeg (D_n)=2n$. Recall that
$$ \GKtrdeg (A):=\sup_{V} \inf_{b} \, \GK (K[bV])$$
where $V$ ranges over the finite dimensional subspaces of $A$ and
$b$ ranges over the regular elements of $A$. Another proofs of the
corollary based on different  ideas were given by A. Joseph
\cite{JosLN74} and R. Resco \cite{Resco79}, see also \cite{MR},
6.6.19. Joseph's proof is based on the fact that the centralizer
of any isomorphic copy of the Weyl algebra $A_n$ in its division
algebra $D_n:=D( \mathbb{A}^n)$ reduces to scalars
(\cite{JosgenQ}, 4.2), Resco proved that $\Cdim (D_n)=n$
(\cite{Resco79}, 4.2) using the result of Rentschler and Gabriel
 \cite{Ren-Gab} that $\Kdim (A_n)=n$  (over an arbitrary field of characteristic
zero).


\section{Filter Dimension  and Isotropic
Subalgebras of Poisson
Algebras}\label{hazisotr}

In this section, we apply Theorem \ref{GKcsuba} to obtain an upper
bound for the Gelfand-Kirillov dimension of  {\bf isotropic}
subalgebras of certain Poisson algebras (Theorem \ref{PoGKcsu}).

Let $(P, \{ \cdot , \cdot \} )$ be a {\bf  Poisson algebra} over
the field $K$. Recall that $P$ is an associative  commutative
$K$-algebra which is a Lie algebra with respect to the bracket $\{
\cdot , \cdot \}$ for which {\em Leibniz's rule} holds:
$$ \{ a, xy\} = \{ a, x\} y +x\{ a, y\}\;\; {\rm for \; all}\;\;
a,x,y\in P,$$
 which means that the {\bf inner derivation } $\ad (a): P\ra P$,
 $x\mapsto \{ a, x\}$, of the Lie algebra $P$ is also a derivation
 of the associative algebra $P$. Therefore, to each Poisson
 algebra $P$ one can attach an associative subalgebra $A(P)$ of
 the ring of differential operators $\CD (P)$ with coefficients
 from the algebra $P$ which is generated by $P$ and $\ad (P):=\{
 \ad (a)\, | \, a\in P\}$. If $P$ is a finitely generated algebra
 then so is the algebra $A(P)$ with $\GK (A(P))\leq \GK (\CD
 (P))<\infty $.

 {\it Example}. Let $P_{2n}=K[x_1, \ldots , x_{2n}]$ be the {\bf
 Poisson polynomial algebra} over a field $K$ of characteristic
 zero equipped with the {\bf Poisson bracket}
 $$ \{ f,g\} =\sum_{i=1}^n( \frac{\der f}{\der x_i}\frac{\der g}{\der
 x_{n+i}}-\frac{\der f}{\der x_{n+i}}\frac{\der g}{\der
 x_i}).$$
 The algebra $A(P_{2n})$ is generated by the elements
 $$x_1, \ldots , x_{2n}, \; \ad (x_i)=\frac{\der }{\der
 x_{n+i}},\;
 \ad (x_{n+i})=-\frac{\der }{\der
 x_{i}}, \; i=1,\ldots , n.$$ So, the algebra $A(P_{2n})$ is canonically
 isomorphic to the Weyl algebra $A_{2n}$.

 {\it Definition}. We say that a Poisson algebra $P$ is a {\bf
 strongly simple Poisson algebra} if
\begin{enumerate}
\item $P$ is a finitely generated (associative) algebra which is a
domain, \item the algebra $A(P)$ is central simple, and \item for
each set of algebraically independent elements $a_1, \ldots , a_m$
of the algebra $P$ such that $\{ a_i, a_j\}=0$ for all $i,j=1,
\ldots , m$ the (commuting) elements $a_1, \ldots , a_m, \ad
(a_1), \ldots ,$ $ \ad (a_m)$ of the algebra $A(P)$ are
algebraically independent.
\end{enumerate}

\begin{theorem}\label{PoGKcsu}
(\cite{dimcom}) Let $P$ be a strongly simple Poisson algebra, and
$C$ be an isotropic subalgebra of $P$, i.e. $\{ C,C\}=0$. Then
$$ \GK (C)\leq \frac{\GK (A(P))}{2}\left( 1-\frac{1}{f_{A(P)}+\max \{
f_{A(P)}, 1\}} \right)$$ where $f_{A(P)}:= \max \{
\fdim_{Q_m}(Q_m\t A(P))\, | \, 0\leq m\leq  \GK (A(P))\}$.
\end{theorem}

{\it Proof}. By the assumption the finitely generated algebra $P$
is a domain, hence the finitely generated algebra $A(P)$ is a
domain (as a subalgebra of the domain $\CD (Q(P))$, the ring of
differential operators with coefficients from the field of
fractions $Q(P)$ for the algebra $P$). It suffices to prove the
inequality for isotropic subalgebras of the Poisson algebra $P$
that are polynomial algebras. So, let $C$ be an isotropic
polynomial subalgebra of $P$ in $m$ variables, say  $a_1, \ldots ,
a_m$. By the assumption, the commuting elements $a_1,\ldots , a_m
, \ad (a_1), \ldots , \ad (a_m)$ of the algebra $A(P)$ are
algebraically independent. So,  the Gelfand-Kirillov dimension of
the subalgebra $C'$ of $A(P)$  generated by these elements is
equal to $2m$. By Theorem \ref{GKcsuba},
$$ 2 \GK (C)=2m=\GK (C')\leq \GK (A(P))\left( 1-\frac{1}{f_{A(P)}+\max \{
f_{A(P)}, 1\}} \right), $$ and this proves the inequality. $\Box $

\begin{corollary}\label{cPoGKcs}
\begin{enumerate}
\item The Poisson polynomial algebra $P_{2n}= K[x_1, \ldots ,
x_{2n}]$ (with the Poisson bracket)  over a field $K$ of
characteristic zero is a strongly simple Poisson algebra, the
algebra $A(P_{2n})$ is canonically isomorphic to the Weyl algebra
$A_{2n}$. \item The Gelfand-Kirillov dimension of every  isotropic
subalgebra of the polynomial Poisson algebra $P_{2n}$ is $\leq n$.
\end{enumerate}
\end{corollary}

{\it Proof}. $1$. The third condition in the definition of
strongly simple Poisson algebra is the only statement we have to
prove. So, let $a_1, \ldots , a_m$ be  algebraically independent
elements of the algebra $P_{2n}$ such that $\{ a_i, a_j\}=0$ for
all $i,j=1, \ldots , m$. One can find polynomials, say
$a_{m+1},\ldots , a_{2n}$,  in $P_{2n}$ such that the elements
$a_1, \ldots , a_{2n}$ are algebraically independent, hence the
determinant $d$ of the Jacobian matrix $J:=(\frac{\der a_i}{\der
x_j })$ is a nonzero polynomial. Let $X=(\{ x_i, x_j\} )$ and
$Y=(\{ a_i, a_j\} )$ be, so-called, the {\em Poisson matrices}
associated with the elements $\{ x_i\}$ and $\{ a_i \} $. It
follows from $Y=J^TXJ$ that $\det (Y)=d^2 \det (X)\neq 0$ since
$\det (X)\neq 0$. The derivations
$$ \d_i:= d^{-1} \det
\begin{pmatrix}
  \{ a_1, a_1\}  & \ldots & \{ a_1 , a_{i-1} \}  & \{ a_1 , \cdot  \}  & \{ a_1 , a_{i+1} \}
   & \ldots  & \{ a_1 , a_{2n} \}  \\
\{ a_2, a_1\}  & \ldots & \{ a_2 , a_{i-1} \}  & \{ a_2 , \cdot \}
& \{ a_2 , a_{i+1} \}
   & \ldots  & \{ a_2 , a_{2n} \}  \\
   &  & & \ldots  &  & &  \\
\{ a_{2n}, a_1\}  & \ldots & \{ a_{2n} , a_{i-1} \}  & \{ a_{2n} ,
\cdot  \} & \{ a_1 , a_{i+1} \}
   & \ldots  & \{ a_{2n} , a_{2n} \}  \\
\end{pmatrix},
$$
$i=1, \ldots , 2n$, of the rational function  field $Q_{2n}=K(x_1,
\ldots , x_{2n})$ satisfy the following properties: $\d_i
(a_j)=\d_{i,j}$, the Kronecker delta. For each $i$ and $j$, the
kernel of  the derivation $\D_{ij}:=\d_i\d_j-\d_j\d_i\in
\Der_K(Q_{2n})$ contains $2n$ algebraically independent elements
 $a_1, \ldots , a_{2n}$. Hence $\D_{ij}=0$ since the field
 $Q_{2n}$ is algebraic over  its subfield $K(a_1, \ldots , a_{2n})$
 and ${\rm char} (K)=0$. So, the subalgebra, say $W$, of the ring
 of differential operators
 $\CD (Q_{2n})$ generated by the elements $a_1, \ldots , a_{2n},
\d_1, \ldots , \d_{2n}$ is isomorphic to the Weyl algebra
$A_{2n}$, and so $\GK (W)=\GK (A_{2n})=4n$.

Let $U$ be the $K$-subalgebra of $\CD (Q_{2n})$ generated by the
elements $x_1, \ldots , x_{2n}, \d_1, \ldots , \d_{2n}$, and
$d^{-1}$. Let $P'$ be the localization of the polynomial algebra
$P_{2n}$ at the powers of the element $d$. Then $\d_1, \ldots ,
\d_{2n}\in \sum_{i=1}^{2n}P'\ad (a_i)$ and $\ad (a_1), \ldots ,
\ad (a_{2n})\in  \sum_{i=1}^{2n}P'\d_i$, hence the algebra $U$ is
generated (over $K$) by $P'$ and  $\ad (a_1) , \ldots , \ad
(a_{2n})$. The algebra $U$ can be viewed as  a subalgebra of the
ring of differential operators $\CD (P')$. Now, the inclusions,
$W\subseteq U\subseteq \CD (P')$ imply $4n=\GK (W)\leq \GK (U)\leq
\GK (\CD (P'))=2\GK (P')=4n$, therefore $\GK (U)=4n$. The algebra
$U$ is a factor algebra of an iterated Ore extension $V=P'[t_1;\ad
(a_1)]\cdots [t_{2n}; \ad (a_{2n})]$. Since $P'$ is a domain, so
is the algebra $V$. The algebra $P'$ is a finitely generated
algebra of Gelfand-Kirillov dimension $2n$, hence $\GK (V)=\GK
(P')+2n =4n$ (by \cite{MR}, 8.2.11). Since $\GK (V)=\GK (U)$ and
any proper factor algebra of $V$ has Gelfand-Kirillov dimension
strictly less than $\GK (V)$ (by \cite{MR}, 8.3.5, since $V$ is a
domain), the algebras $V$ and $U$ must be isomorphic. Therefore,
 the (commuting) elements $a_1, \ldots
, a_m, \ad (a_1), \ldots ,$ $ \ad (a_m)$ of the algebra $U$ (and
$A(P)$) must be  algebraically independent.

$2$. Let $C$ be an isotropic subalgebra of the Poisson algebra
$P_{2n}$. Note that $f_{A(P_{2n})}=f_{A_{2n}}=1$  and $\GK
(A_{2n})=4n$. By Theorem \ref{PoGKcsu},
$$ \GK (C) \leq \frac{4n}{2}(1-\frac{1}{1+1})=n. \; \Box $$

{\it Remark}. This result means that for the Poisson polynomial
algebra $P_{2n}$  the right hand side of the inequality of Theorem
\ref{PoGKcsu} is the {\em exact} upper bound for the
Gelfand-Kirillov dimension of isotropic subalgebras in $P_{2n}$
since the polynomial subalgebra $K[x_1, \ldots, x_n]$ of $P_{2n}$
is isotropic.

Department of Pure Mathematics

University of  Sheffield

Hicks Building

Sheffield S3~7RH

UK

email: v.bavula@sheffield.ac.uk



\end{document}